\documentclass[12pt]{article}
\usepackage[usenames]{color}
\newcounter{biblio}
\newenvironment{references}%
{\begin{list}{[\arabic{biblio}]}{\usecounter{biblio}%
\setlength{\leftmargin}{2.5em}\setlength{\rightmargin}{0pt}%
\setlength{\labelwidth}{2em}\setlength{\itemsep}{0pt}}}{\end{list}}
\newcommand{\References}%
{\vspace{2.8ex plus .3ex minus .3ex}%
\begin{center}{\bf References}\end{center}\begin{references}}

\usepackage{latexsym}
\usepackage{graphicx}
\usepackage{amssymb, amsmath}
\usepackage{mathrsfs}
\newcommand{\N}{{\mathbb{N}}}
\newcommand{\Z}{{\mathbb{Z}}}

\newcommand{\zd}{\Z^d}

\newcommand{\R}{{\mathbb{R}}}

\newcommand{\vvs}{\vspace{2ex}}

\newcommand{\lef}{\left}
\newcommand{\rig}{\right}

\newcommand{\8}{\infty}

\renewcommand{\b}{\beta}

\newcommand{\h}{\eta}

\newcommand{\lm}{\lambda}

\newcommand{\n}{\nu}

\newcommand{\s}{\sigma}

\renewcommand{\t}{\tau}

\newcommand{\eps}{\epsilon}

\newcommand{\cF }{{\cal F}}
\newcommand{\cG }{{\cal G}}

\newcommand{\cN }{{\cal N}}

\newcommand{\cS }{{\cal S}}

\makeatletter
\def\section{\@startsection{section}{1}{\z@}{-3.5ex plus -1ex minus
 -.2ex}{2.3ex plus .2ex}{\bf}}
\def\subsection{\@startsection{subsection}{2}{\z@}{-3.25ex plus -1ex minus
 -.2ex}{1.5ex plus .2ex}{\bf}}
\makeatother

\newtheorem{theorem}{Theorem}[section]
\newtheorem{lemma}[theorem]{Lemma}
\newtheorem{proposition}[theorem]{Proposition}
\newtheorem{corollary}[theorem]{Corollary}

\newtheorem{remark}[theorem]{Remark}

\newcommand{\tS}{\widetilde S}
\newcommand{\E}{{\mathbb E}}
\newcommand{\Var}{{\rm Var}}
\newcommand{\Cov}{{\rm Cov}}

\newcommand{\Cstval}{\frac{d \;\Gamma(d/2)}{(d-2)\pi^{d/2}}}
\newcommand{\Cst}{{\mathfrak K}_d}
\newcommand{\Zd}{{\mathfrak Z}_d}
\newcommand{\Cstd}{{\kappa_2}}
\newcommand{\CstF}{{\mathfrak C}_d}

\newcommand{\cvlaw}{\stackrel{\rm{ law}}{\longrightarrow}}

\newcommand{\nn}{\nonumber}

\def\qed{\hfill\rule{.2cm}{.2cm}\par\medskip\par\relax}

\begin{document}

\vvs

\title{Rate of convergence for polymers in a weak disorder}
\author{Francis Comets$^{1}$ \and
 Quansheng Liu$^{2}$}

\maketitle

{\footnotesize
\noindent $^{~1}$Universit\'e Paris Diderot -- Paris 7,
Math\'ematiques,
 case 7012, F--75205 Paris
Cedex 13, France
\\
\noindent e-mail:
\texttt{comets@math.univ-paris-diderot.fr}

\noindent $^{~2}$  Univ. Bretagne-Sud, UMR 6205, LMBA, F-56000 Vannes, France\\
\noindent e-mail: \texttt{quansheng.liu@univ-ubs.fr}
}

\begin{abstract}
We consider directed polymers in random environment on the lattice $\Z^d$ at small inverse temperature and dimension $d \geq 3$.
Then, the normalized partition function $W_n$ is a regular martingale with limit $W$.
We prove that $n^{(d-2)/4} (W_n-W)/W_n$ converges in distribution to
a Gaussian law. Both the  polynomial rate of convergence and the scaling with the  martingale $W_n$ are different from those for polymers on trees.
 \\[.3cm]\textbf{Keywords:} Directed polymers, random environment, weak disorder, rate of convergence, regular martingale, central limit theorem for martingales, stable and mixing convergence.
\\[.3cm]\textbf{AMS 2010 subject classifications:}
Primary 60K37. Secondary 60F05, 60G45, 60J80, 82D60.
\end{abstract}

%
%
%
%

\section{Polymer  models and statement of the main result}

\subsection{Motivation}

We consider directed polymers  in random environment, given by a simple random walk on the
$d$-dimensional lattice in a space-time random potential.
In a seminal paper, Derrida and Spohn~\cite{DerridaSpohn88} perform a detailed  analysis of
polymers on the Cayley tree, or equivalently, the branching random walk with a fixed branching number.
Later the same model has been taken up as an approximation and a toy model with explicit computations:  in  the physics literature, we mention the pleasant, recent and documented
survey~\cite{HalpinHealyTakeuchi15},
and also \cite{DeanMajumdar01} for the statistics of extremes on the hierarchical tree at zero temperature;  on the mathematical side, the authors of~\cite{AlbertsOrtgiese13} study the near-critical scaling window on the tree, the analogue of the intermediate disorder regime where the rescaled lattice model on line converges to the KPZ continuum random polymer \cite{AKQ14, CaSuZyJEMS}.
Not only a source of inspiration and guidance,
this model, as well as related random cascades,
were also found  to provide quantitative bounds on polymer models on the lattice in~\cite{CVargas06,LiuWatbled, Nguyen16}.

In spite of these similarities, the two models
behave quite differently in many  aspects. In the strong disorder phase, the free energy
of the branching process is linear in the inverse temperature $\beta$ though it is strictly convex for
the polymer on the lattice, see Theorem 1.5 in~\cite{CPopovVachkovskaia08} in the case of a Bernoulli
environment. Also, the fluctuations are expected to be of a completely different nature in the two models. In this paper,
we consider the weak disorder regime, and we show that the martingale convergence takes place at a polynomial rate, whereas it is exponential
in the corresponding supercritical  Galton-Watson process \cite{Heyde, Heyde71}.

More precisely, it is shown in \cite{Heyde, Heyde71} that, for  a Galton-Watson process $(Z_n)$  with $Z_0=1$,
$m= E Z_1 >1$ and $E Z_1^2 < \infty$, the renormalized population size $W_n= Z_n/ m^n$ is a regular martingale with limit $W$   such that
\begin{equation}
\label{eq:GW}
m^{n/2} (W-W_n) \to a W^{1/2} G \quad \mbox{ in distribution }
\end{equation}
and
\begin{equation}
\label{eq:GW2}
 m^{n/2} \, \frac{(W-W_n)}{ W_n^{1/2}}  \to a G \quad \mbox{ in distribution},
\end{equation}
where $a^2 = \frac{\Var Z_1}{m^2-m}$,
$G  $  is a Gaussian $ \cN (0,1 )$ distributed random variable independent of   $W$.
Similarly, for branching random walks, the convergence of the Biggins martingale to its limit is exponentially fast \cite{Iksanov06, IksanovKabluchko15} in the regular case.
Recently the same question was studied   for a branching process in a random environment \cite{HuangLiu14a,WangGaoLiu11},
leading  to similar conclusions.

In this paper, we consider  random polymers on the lattice in a time-space dependent random medium, deep inside the weak disorder regime.
Similar to the supercritical case of a branching process, weak disorder can be defined as the regime where the natural martingale is regular \cite{Bol89, ImSp88}, or where
the polymer is diffusive \cite{CoYo06}. It holds in space dimension  $d\geq 3$ \cite{Lacoin10} and at a temperature larger than some critical value which can be estimated
by second moment and entropy considerations \cite{Bir04, CaCa09, HuSh09}.
In Theorem \ref{th:clt} below, we prove that,
at large temperature, the speed of convergence is polynomial but not exponential, and the limit scales with $W$ or $W_n$  instead of their square root
as in (\ref{eq:GW}) and (\ref{eq:GW2}). Precisely, we show a central limit theorem for the difference between the martingale and its limit:  the ratio of the difference
divided by $n^{-(d-2)/4}$ times the martingale is asymptotically normal.

 In view of  (\ref{eq:GW}) and (\ref{eq:GW2}), this limit has two remarkable and unexpected features.
The slowdown in the rate of convergence (compared to the branching case) is due to space correlations coming from further intersections between paths on the lattice
but not on the tree. Also the unusual linear scaling in the martingale can be understood as coming from fluctuations, and quadratic variations scale like the square of the martingale.

\subsection{Notations}

$\bullet$ {\it The random walk:}
$(\{ S_n\}_{n \geq 0}, P_x)$ is a nearest neighbor, symmetric simple random walk on
the $d$-dimensional integer lattice $\Z^d$ starting from $x$, $d \geq 3$.
We let $P=P_0$  and we denote by $P[ f ] = \int f \, dP $ the expectation of $f$ with respect to $P$.
\\
$\bullet$ {\it The random environment:}
$\h  =\{\h (n,x) : n \in \N,\; x \in \zd \}$ is an  independent and identically distributed (i.i.d.)
sequence of real random variables (r.v.'s),
non-constant, such that,
\begin{equation} \label{expint} \nn
 \lm(\b)  :=
 \ln \E[\exp (\b \h (n,x))] <\8 \; \; \; \mbox{for all $\b \in \R$,}
\end{equation}
where we denote by $\E$ the  expectation over the environment. The corresponding probability 
measure will be denoted by $\mathbb P$. 
\\
$\bullet$ The {\it partition function} at  inverse temperature $\b \in \R$
\begin{equation} \label{Wn}
W_n
:=P\lef[\exp \lef( \beta H_n(S) - n \lm (\beta) \rig) \rig] \; \mbox{ with } \; H_n(S)= \sum_{0 \leq t \leq n-1} \h (t,S_t),   \quad n\geq 1,
\end{equation}
is the  normalization which makes the Gibbs measure
${W_n}^{-1}
\exp \{\beta H_n(S)-n \lm (\beta) \} \;dP$ a probability measure on the path space for a fixed realization of the environment.
Note that we use a slightly different definition than usual, including $t=0$ but not $t=n$ in the Hamiltonian $H_n$.
This makes no fundamental difference in the results (see Remark \ref{rk:pasdediff} below), but it yields simpler formulas here: 
for two independent simple random walks $S=(S_t)$ and $\tS=(\tS_t)$, we have 
${\rm Cov}_{\mathbb P}(H_n(S), H_n(\tS))=\Var(\eta(0,0) ) \,  N_n$,  where ${\rm Cov}_{\mathbb P}$ denotes the covariance with respect to 
   $\mathbb P$, $N_n$ is  the number of intersections of the paths
$S, {\tilde S}$ up to time $n$:
\begin{equation} \label{Nn}
N_n=N_n(S, {\tilde S}):= \sum_{t=0}^{n-1}  {\bf 1}_{S_t={\tilde S}_t},
\end{equation}
whose limit
\begin{equation}\label{def-Ninfty}
  N_\infty = N_\infty (S, {\tilde S}) := \sum_{t=0}^{\8}  {\bf 1}_{S_t={\tilde S}_t}
\end{equation}
has expectation given by the standard Green function \eqref{def:Green}.
The sequence $(W_n)$ depends on the environment, and it is a positive
martingale with respect to the filtration 
\begin{equation} \label{def-algGn}
\cG_n=\sigma \{\h(t,x); t \leq n-1, x \in \Z^d\}, \quad n\geq 1.
\end{equation}
It is well known  \cite{Bol89, ImSp88}  that
$$
W=\lim_{n \to \8} W_n \quad {\rm exists \ a.s.,\ with } \quad
{\mathbb P}(W>0 )=0 {\rm \ or \ }1.
$$
Moreover, the convergence holds in $L^2$ for $\b$ in a neighborhood of $0$, defined by
\begin{equation} \label{L2}
{\rm {\bf (L2)}}\qquad
 \lambda_2 := \lm (2 \b) - 2\lm ( \b) < \ln (1/\pi_d),
\end{equation}
where  $\pi_d$ is the return probability  of the simple random walk,
\begin{equation} \label{rec/tra}
\pi_d :=
P\{ \mbox{$S_n = 0$ for some $n \geq 1$}\} \in (0,1)
\end{equation}
by transience since $d \geq 3$.
Now, we give a short account of the main steps of the computation of \cite{Bol89}, which is useful for the sequel. 
We can express $W_n^2$ as a sum $W_n^2=
 P^{\otimes 2} \left[ e^{\b [H_n(S)+
H_n({\tilde S})]-2n \lm(\b)}\right]$ over independent paths (so-called replicas), and
we compute, using Fubini's theorem  and independence,
\begin{eqnarray} \nn
\E [W_n^2]
&=&
P^{\otimes 2} \left[ \prod_{t=0}^{n-1}  \E\; e^{\b [\h(t,S_t)+
\h(t,{\tilde S}_t)]-2 \lm(\b)}\right]\\ \nn
&=&  P^{\otimes 2}\left[ \prod_{t=0}^{n-1} \left(e^{\lm(2 \b)-2 \lm(\b)}
 {\bf 1}_{S_t={\tilde S}_t} +  {\bf 1}_{S_t \neq {\tilde S}_t}\right)\right]
\\ \nn
&=&  P^{\otimes 2}\left[ \prod_{t=0}^{n-1} e^{\lm_2
 {\bf 1}_{S_t={\tilde S}_t}} \right] \\ \label{eq:2ndmoment}
&=&  P^{\otimes 2} \left[ e^{ \lambda_2 N_n} \right]\;, \label{eq:NYC}
\end{eqnarray}
with $N_n$ as in \eqref{Nn}.
As $n \to \8$, $N_n \nearrow N_\8  = \sum_{t=0}^{\infty}  {\bf 1}_{S_t={\tilde S}_t} $.
Since the process $(S_t-{\tilde S})$ under $ P^{\otimes 2} $ has the same law as $(S_{2t})$ under $P$ and
$S_{2t+1} \neq 0 \; P$ -a.s.,
$N_n$ has the same law as the number of visits to $0$ of a simple random walk in time $2n$, and
 $N_\8$   is
geometrically distributed with "failure" probability $\pi_d$:
\begin{equation}\label{eq-distri-Ninfty}
   P^{\otimes 2} ( N_\8 = k) = \pi_d^{k-1} (1-\pi_d) \quad \mbox{ for }  k \geq 1.
\end{equation}
 Therefore, by the monotone convergence theorem, as $n \to \8$,
\begin{eqnarray*}
\E[W_n^2] \nearrow P^{\otimes 2} \left[ e^{ \lambda_2 N_\8}\right]
&=& \left\{
\begin{array}{ll}
\frac{\displaystyle (1-\pi_d)e^{\lambda_2}}{\displaystyle 1-\pi_d e^{\lambda_2}} & {\rm if } \; \pi_d e^{\lambda_2}<1, \\
+\8  & {\rm if } \; \pi_d e^{\lambda_2} \geq 1.
\end{array}
\right.
\end{eqnarray*}
 In this paper, we always assume (\ref{L2}), so that
 the martingale  $(W_n)$  is
bounded in $L^2$.
By Doob's convergence theorem, $W_n \to W$ in
$L^2$ and then $W>0$. In particular,
\begin{equation}
\label{eq:varW}
\E W^2 = \frac{\displaystyle (1-\pi_d)e^{\lambda_2}}{\displaystyle 1-\pi_d e^{\lambda_2}}
 \;\; \mbox{ and } \;\;
 \Var \, (W)
 = \frac{\displaystyle  e^{\lambda_2}-1}{\displaystyle 1-\pi_d e^{\lambda_2}} .
\end{equation}

\subsection{A Gaussian limit and the rate of convergence}


Before coming to our main result, we recall two convergence modes. Let $(Y_n)$ be a sequence of real random variables defined on a common probability space
$(\Omega, {\cal F}, P)$, converging in distribution to a limit $Y$.
\begin{itemize}
\item This convergence is called {\em stable} if for all $B \in \cal F$ with $P(B)>0$, the conditional law of $Y_n$ given $B$ converges to some probability distribution
depending on $B$.
\item This convergence is called {\em mixing} if it is stable and the limit of conditional laws does not depend on $B$ -- and therefore is the law of $Y$ --.
\end{itemize}
The stable convergence allows to add extra variables:
for
any fixed r.v. $Z$ on $(\Omega, {\cal F}, P)$, the couple $(Y_n,Z)$ converges in law to some coupling of
$Y$ and $Z$ on an extended space.
The mixing convergence means that $Y_n$ is asymptotically independent of all event $A  \in {\cal F}$. These convergences were introduced by R\'enyi~\cite{Renyi63};
we refer to \cite{AldousEagleson78} for a nice presentation
with the main consequences, and to  \cite{Hall-Heyde} pp. 56-57 for  an extended account on the connections with martingale central limit theorem.

\begin{theorem} \label{th:clt} For $d \geq 3$, there exists some $\beta_0>0$ such that, for $|\beta| < \beta_0$,
\begin{equation} \label{clt2bis1}
 n^{ \frac{d-2}{4}}    (W-W_n)  \rightarrow  \sigma_1 W G
     \quad \mbox{ in distribution}
\end{equation}
and
\begin{equation} \label{clt2bis}
n^{ \frac{d-2}{4}}
 \, \frac{  (W-W_n)}{ W_n} \rightarrow \sigma_1 G
     \quad \mbox{ in distribution, }
\end{equation}
where  $\sigma_1$ is from (\ref{eq:sigma1}),
 $G$ is  a Gaussian r.v. with law $N(0, 1)$, which is
 independent of $W$.
 Moreover, the convergence in \eqref{clt2bis1} is stable, and the convergence in \eqref{clt2bis} is  mixing.

\end{theorem}
The theorem calls for some comments. The value of $\beta_0$ is defined by the conditions in Lemma \ref{lem:vargas}, Lemma \ref{lem:moments} (b) and
Lemma \ref{lem:rate-conv-Dk}.
The result is quite different from \eqref{eq:GW}--\eqref{eq:GW2}.
The speed of convergence of $W_n$ to its limit is $ n^{ -\frac{d-2}{4}} $.
The mixing convergence in \eqref{clt2bis}
shows that the random variable $G_n$ defined in the left hand side of $ \eqref{clt2bis} $ is asymptotically independent of each event $A$ of the environmental probability space, in the sense
that for all real $y$, we have $ \lim_{n\to \8} P (\{G_n \leq y \} \cap A ) = P (G \leq y) P (A)$.
Our approach relies on a central limit theorem for martingales.
\begin{remark}[Usual Hamiltonian] \label{rk:pasdediff}
For the standard Hamiltonian collecting the environment at times $1,2,\ldots n$,  
$$
\overline{W}_n
:=P\lef[\exp \lef(  \beta \sum_{1 \leq t \leq n} \h (t,S_t) - n \lm (\beta) \rig) \rig] = W_{n+1} \exp\{-\beta \h(0,0) + \lambda(\beta)\}\;,
$$
it is straightforward to see that, for $d\geq 3$ and $|\beta| < \beta_0$,  $\overline{W}_n \to \overline{W}:= W \exp\{-\beta \h(0,0) + \lambda(\beta)\}$, that
\begin{equation} \label{clt2bis1t} \nn
 n^{ \frac{d-2}{4}}    ( \overline{W}- \overline{W}_n)  \rightarrow  \sigma_1  \overline{W} G
     \quad \mbox{ in distribution}
\end{equation}
and
\begin{equation} \label{clt2bist} \nn
n^{ \frac{d-2}{4}}
 \, \frac{  ( \overline{W}- \overline{W}_n)}{  \overline{W}_n} \rightarrow \sigma_1 G
     \quad \mbox{ in distribution, }
\end{equation}
with  $\sigma_1$ and
 $G$ as above and $G$ independent from $\overline{W}$.
\end{remark}

\medskip

{\bf Organization of the paper:} In Section \ref{sec:2} we start with algebraic computations of covariances and use Green function estimates to derive asymptotics. Of independent interest, a specific form of central limit theorem for infinite martingales arrays is given in Lemma \ref{lem:martingale-clt} of Section \ref{sec:3}, and used to prove Theorem \ref{th:clt}.
The proofs of intermediate lemmas are postponed to Section \ref{sec:4}. A key step consists in controling the sum of conditional variances in the quadratic norm, so we implement the necessary algebra for a system of 4 replicas at the beginning of the section.


\section{The correlation structure} \label{sec:2}
It is useful to introduce
\begin{equation}\label{eq:Wn(x)}
W_n(x)= P \left[ e^{\b H_n(S)
-n \lm(\b)} {\bf 1}_{S_n=x}\right],
\end{equation}
and to observe that, for $m \geq 0$ including $m=\8$ with the convention $W_\8=W$,
\begin{equation}
\label{eq:decMov}
W_{n+m}= \sum_x W_n(x) W_m \circ \theta_{n,x},
\end{equation}
by  Markov property. We view $W_n$ as a function of $\eta$, we denote by $\theta_{n,x}$  the shift operator on the environment,
$\theta_{n,x}\eta: (t,y) \mapsto \eta(n+t,x+y)$. We use $\theta_x=\theta_{0,x}$ as a short notation.
Taking $m=\8$ we see that
\begin{equation} \label{eq:diff}
W-W_n=  \sum_x W_n(x) \big( W \circ \theta_{n,x}-1\big) .
\end{equation}
For nearest-neighbor paths $S, \tS$ define $\tau_n$ the time delay of first intersection after time
$n \geq 0$,
\begin{equation}\nn
\tau_n(S,\tS)=\inf\{k \geq 0: S_{n+k}=\tS_{n+k}\} ,
\end{equation}
with the convention that $\inf \emptyset = +\8$. We write $\t=\t_0$, and we note that $\pi_d=
P_{0,0}^{\otimes 2} ( \t_1 <\8)$.

\subsection{Covariance of the martingale limit and rate of convergence in $L^2$ }

Recall that $ \Var(W)  =  \frac{\displaystyle  e^{\lambda_2}-1}{\displaystyle 1-
P_{0,0}^{\otimes 2} ( \t_1 <\8)
 e^{\lambda_2}}$ , cf. \eqref{eq:varW}.
\begin{proposition} \label{prop:2.1}
Assume $\lambda_2 < \frac{2}{d} \ln \frac{1}{\pi_d}$. Then,
\begin{equation} \label{eq:covW}
\Cov (W, W\circ \theta_x)= \Var(W)   \times P_{0,x}^{\otimes 2} (\tau <\8)\;.
\end{equation}
Moreover,
\begin{eqnarray} \label{eq:dist2}
\|W-W_n\|_2^2 &=& \Var (W) \times  P_{0,0}^{\otimes 2} (e^{\lambda_2 N_n} {\bf 1}_{\t_n<\8}) ,
\end{eqnarray}
and, as $n \to \8$ ,
\begin{eqnarray}
\label{eq:dist2cv}
 \|W-W_n\|_2^2
&\sim &
\sigma^2 \times n^{-(d-2)/2} , 
\end{eqnarray}
with the constant $\s^2$ from (\ref{eq:sigma}).
\end{proposition}
{\bf Proof.}
We first compute the covariance of $W$ and $W\circ \theta_x$. Denote by $\cF_n$ the $\sigma$-field generated by $S_i, \tS_i$ for $0\leq i \leq n$. By convergence in $L^2$,
we obtain as in \eqref{eq:NYC},
\begin{eqnarray} \nn
\Cov (W, W\circ \theta_x)&=&
\lim_{m \to \8} \E\big[(W_m-1)(W_m\circ\theta_{x}-1)\big]  \qquad 
\\ \nn
&=&
\lim_{m \to \8} P_{0,x}^{\otimes 2}
\E\big[(e^{\b H_m(S)-m \lm(\b)}-1)(e^{\b H_m(\tS)-m \lm(\b)}-1)\big] \qquad {\rm (Fubini)}\\ \nn
&=&
\lim_{m \to \8} P_{0,x}^{\otimes 2}
\E\big[e^{\b H_m(S)+\b H_m(\tS)-2m \lm(\b)}-1\big] \\ \nn
&=&
\lim_{m \to \8} P_{0,x}^{\otimes 2}
\big[e^{\lambda_2 N_m}\big] -1\\ \nn
&=&
 P_{0,x}^{\otimes 2}
\big[(e^{\lambda_2 N_\8}-1)\big] \\ \nn
&=&
 P_{0,x}^{\otimes 2} \Big[ {\bf 1}_{\tau < \8}
 P_{0,x}^{\otimes 2}
\big[(e^{\lambda_2 N_\8}-1) \vert {\cF_\tau} \big] \Big]\
\\ \nn 
&=& P_{0,x}^{\otimes 2} (\tau <\8)
 P_{0,0}^{\otimes 2}
\big[(e^{\lambda_2 N_\8}-1)\big]  \qquad\qquad {\rm (strong\ Markov\  property)}
\end{eqnarray}
which is the first claim since $\Var(W)=P_{0,0}^{\otimes 2}
\big[e^{\lambda_2 N_\8}\big] -1$.

We next calculate the $L^2$ norm of $W-W_n$.
By (\ref{eq:diff}),
\begin{eqnarray*}
\|W-W_n\|_2^2   
&=& \E \Big(\sum_x W_n(x) (W\circ \theta_{n,x}-1)\Big)^2  \\
&=& \E \Big(\sum_{x,y} W_n(x) W_n(y) (W\circ \theta_{n,x}-1) (W\circ \theta_{n,y}-1)\Big)\\
&=& \sum_{x,y} \E [W_n(x) W_n(y)] \E[(W\circ \theta_{n,x}-1) (W\circ \theta_{n,y}-1)]
 \qquad({\rm  independence}) \\
&=& \Var(W)   \sum_{x,y} \E \left[ W_n(x) W_n(y)\right]
 P_{x,y}^{\otimes 2} (\tau <\8) \qquad({\rm by\ } (\ref{eq:covW}))\\
 &=& \Var(W)   \sum_{x,y} P_{0,0}^{\otimes 2}  \left[ e^{\lambda_2 N_n} {\bf 1}_{ \{S_n=x, \tS_n=y\}}\right]
 P_{x,y}^{\otimes 2} (\tau <\8)\qquad  ({\rm cf.\ } \eqref{eq:NYC})
  \\
 &=& \Var (W) \times  P_{0,0}^{\otimes 2} (e^{\lambda_2 N_n} {\bf 1}_{\t_n<\8}),
\end{eqnarray*}
by Markov property. This is (\ref{eq:dist2}).

We finally derive  (\ref{eq:dist2cv}) using classical estimates for the Green function.
Denote by
\begin{equation} \label{def:Green}
G(x) : =P_{0,x}^{\otimes 2}\big[ \sum_{n=0}^\8 {\bf 1}_{\{S_n-\tS_n=0\}}\big] = P_{0,x}^{\otimes 2} \; (N_ \8)
\end{equation}
 the Green function for the symmetrized walk $(S_n-\tS_n)_n$, and observe that it is equal on even sites $x$ (i.e., when $\|x\|_1=0$ mod 2) to the Green function of the simple random walk:
   for even sites $x$
  $$ G(x) = P_{x}\big[ \sum_{n=0}^\8 {\bf 1}_{S_n=0}\big],$$
  since the process $ (S_n-\tS_n)$ under $P_{0,x}^{\otimes 2} $ has the same law as $(S_{2n})$ under $P_{x}$ and
  $S_{2n+1}\neq 0$ $P_{x}$ a.s. (for even sites $x$). For odd sites $x$, $G(x) =0$, since in this case $ S_n-\tS_n \neq 0 $ $P_{0,x}^{\otimes 2} $ a.s..
   From the geometric distribution of $N_\8$ under $ P_{0,0}^{\otimes 2}$,   we have
 $G(0)=(1-\pi_d)^{-1}$. By Markov property we have
 $$
 P_{0,x}^{\otimes 2}(\tau <\8)= G(x)/G(0), \qquad x \in \Z^d,
 $$
 and can use the  classical estimates for the Green function, see e.g. \cite[Theorem 4.3.1]{LawlerLimic}:
  for even sites $x$,
 \begin{equation}
\label{eq:greenest}
G(x) = \frac{\Cst}{|x|^{d-2}} + {\mathcal O}\big(\frac{1}{|x|^d}\big), \qquad |x| \to \8,
\end{equation}
where $\Cst \in (0,\infty) $ is a constant whose   value is
\begin{equation}
\label{eq:valeurconstante}
\Cst=\Cstval.
\end{equation}
By the central limit theorem, we have the
  following
 convergence in distribution under $P_{0,0}^{\otimes 2}$  to a Gaussian vector:
$$
n^{-1/2} (S_n-\tS_n) \to Z \quad \mbox{ in distribution, }
$$
where $Z$ is a Gaussian vector with mean $0$  and covariance matrix $(2/d)I_d$.
 Thus, with (\ref{eq:greenest}),
\begin{eqnarray} \label{eq:NYC2}
n^{(d-2)/2} G(S_n-\tS_n) \to \frac{\Cst}{|Z|^{d-2}} \quad \mbox{ in distribution.}
\end{eqnarray}

We shall need the following two lemmas, whose proofs are postponed by the end of the section.

\begin{lemma}[Asymptotic independence]
\label{lem:cvloi} 
Under $P^{\otimes 2}$ we have
the following  convergence in distribution  of random vectors:
$$
\Big(N_n, n^{-1/2} S_n, n^{-1/2}\tS_n\Big) \longrightarrow (N, Z_1,\widetilde Z_1) \quad \mbox{ in distribution, }
$$
where $N,Z_1,\widetilde Z_1$ are independent, and where
\begin{itemize}
\item $N$ is geometrically distributed with parameter $1-\pi_d$:
$P(N=k)=\pi_d^{k-1} (1-\pi_d)$ for $k \geq 1$,
\item  $Z_1$ and $\widetilde Z_1$ are Gaussian vectors
with mean $0$  and covariance matrix $(1/d)I_d$.
\end{itemize}
\end{lemma}

\begin{lemma}[Boundedness in $L^{1+\delta}$]
\label{lem:moments} (a) For $a>0$,
  \begin{equation} \label{eqn:moments}
  \limsup_{n\rightarrow \infty} P^{\otimes 2}
   \big( \Big| \frac{ S_n -\tilde S_n}{\sqrt n}  \Big|^{-a}
    {\bf 1}_{\{  S_n - \tilde S_n \neq 0 \}}\big) <\infty \quad \mbox{ if and only if } a<d.
    \end{equation}
    (b) If $\lambda_2 < \frac{2}{d} \ln \frac{1}{\pi_d}$, then for $\delta >0$ small enough,
    \begin{equation} \label{eqn:momentsG}
      \limsup_{n\rightarrow \infty} P^{\otimes 2} [  e^{\lambda_2 N_n }
      n^{\frac{d-2}{2}} G(S_n - \tilde S_n) ]^{1+\delta} < \infty.
   \end{equation}
\end{lemma}

We first end the proof of Proposition \ref{prop:2.1}.
From Lemma \ref{lem:cvloi}, with $Z=Z_1-\widetilde Z_1$,
\begin{equation}
\label{eq:cvloi11}
e^{\lambda_2 N_n}\times  n^\frac{d-2}{2} G(S_n-\tS_n) \cvlaw
e^{\lambda_2 N}\times   \frac{\Cst}{|Z|^{d-2}} .
\end{equation}
By Lemma \ref{lem:moments}(b),
 the sequence in the left-hand side of (\ref{eq:cvloi11})
is uniformly integrable, so that the  convergence in law implies the convergence of expectations, allowing to write
 \begin{eqnarray} \nn
n^\frac{d-2}{2}P_{0,0}^{\otimes 2} (e^{\lambda_2 N_n} {\bf 1}_{\t_n<\8}) &=&
n^\frac{d-2}{2}P_{0,0}^{\otimes 2} \big[e^{\lambda_2 N_n} P_{0,0}^{\otimes 2} ( \t_n<\8 | \cF_n) \big]\\ \nn
 &=&
G(0)^{-1} P_{0,0}^{\otimes 2} \left[e^{\lambda_2 N_n}\times  n^\frac{d-2}{2} G(S_n-\tS_n) \right]
\\ \label{eq:grosqnu}
 &\longrightarrow&
G(0)^{-1} P_{0,0}^{\otimes 2} \left[e^{\lambda_2 N_\8}\right]
\times  E\left[     \frac{\Cst}{|Z|^{d-2}}  \right] .
\end{eqnarray}
Together with \eqref{eq:dist2},
this ends the proof of \eqref{eq:dist2cv}, yielding the value
\begin{eqnarray}
\label{eq:sigma2} \nn
  \sigma^2&=&  \frac{\Cst \Zd}{G(0)} \Var(W)  \E(W^2)\\&=&   {\Cst \Zd (1-\pi_d)^2 } \times
 (e^{\lambda_2}-1)e^{\lambda_2} \times \frac{1}{\left(1-\pi_d e^{\lambda_2}\right)^{2}}, \label{eq:sigma}
\end{eqnarray}
 with $\Var (W)$ from (\ref{eq:varW}) and
 \begin{equation} \label{val:Zd}
	\Zd \stackrel{\rm def}{=}  E\left[   \frac{1}{|Z|^{d-2}}  \right] = \frac{1}{\Gamma(d/2)}  \left(\frac{d}{4}\right)^{(d-2)/2}  
\end{equation}
from the chi-square distribution.
This proof of
Proposition \ref{prop:2.1} is complete.
\qed
For later purposes, define
\begin{equation}
\label{eq:sigma1}
\sigma_1^2= \Cst \Zd (1-\pi_d) \times \Var (W),
\end{equation}
so that $\sigma^2= \sigma_1^2 \E W^2$.

\subsection{Proof of the lemmas}
It remains to prove Lemmas \ref{lem:cvloi} and \ref{lem:moments}.

\noindent
{\bf Proof of Lemma \ref{lem:cvloi}}.
First observe  that
$$\sup_{n \geq m} P^{\otimes 2}(N_n > N_m) \to 0 \;\; {\rm as } \; m \to \8,$$
since $N_n \nearrow N_\8 < \8$ a.s.
Fix $m \geq 1$ and  $f,g,\tilde g$ continuous and bounded. For all $n \geq m$,
we write
\begin{eqnarray*}
\!\!\!\!\!\!\!\!\!\!\!\!\!\!\!\!\!\!\!\!\!
&&
\!\!\!\!\!\!\!\!\!\!\!\!\!\!\!\!\!\!\!\!\!
P^{\otimes 2}   \left[ f(N_n) g\big( n^{-1/2}S_n
\big) \tilde
  g\big( n^{-1/2}\tilde S_n
\big)
\right]
\qquad \qquad \qquad \qquad \qquad \qquad \qquad
 \\
&=&
P^{\otimes 2}  \left[ f(N_n) g\big( n^{-1/2}S_n
\big) \tilde
  g\big( n^{-1/2}\tilde S_n
\big)  {\bf 1}_{N_n = N_m}
\right]
 + \eps(n,m)\\
&=&
P^{\otimes 2}  \left[ f(N_m) g\big( n^{-1/2}S_n
\big) \tilde
  g\big( n^{-1/2}\tilde S_n
\big)  {\bf 1}_{N_n = N_m}
\right]
 + \eps(n,m)
\\
&=&
 P^{\otimes 2}  \left[ f(N_m) g\big( n^{-1/2}(S_n-S_m)
\big) \tilde
  g\big( n^{-1/2}(\tilde S_n-\tilde S_m)
\big)  {\bf 1}_{N_n = N_m}
\right]
 + \eps'(n,m)\\
&=&
 P^{\otimes 2} \left[ f(N_m) g\big(  n^{-1/2}(S_n-S_m)
\big) \tilde
  g\big(  n^{-1/2}(\tilde S_n-\tilde S_m)
\big)
\right]
+ \eps"(n,m)
\\
&=&
 P^{\otimes 2} [ f(N_m)] \times P\left[ g\big(  n^{-1/2}(S_n-S_m)
\big)\right] \times
 P\left[  \tilde g\big(  n^{-1/2}(\tilde S_n-\tilde S_m)
\big)
\right]
+ \eps"(n,m)\;,
\end{eqnarray*}
where the equalities define the terms $\eps(n,m), \eps'(n,m), \eps''(n,m)$
on their first occurrence.
Here,
$$|\eps(n,m)| \leq \|f\|_\8  \|g\|_\8  \|\tilde g \|_\8
P(N_n \neq N_m)$$ tends to 0 as $m \to \8$ uniformly in $n \geq m$,
$\eps'(n,m)-\eps(n,m) \to 0$ as $n \to \8$ for all fixed $m$,
and $\sup_{n \geq m}\eps''(n,m) \to 0$ as $m \to \8$. The last equality
comes from independence in the increments of the random walks, and of the
two random walks $S$ and $\tilde S$. Hence, letting $n \to \8$ and then
$m \to \8$, we get
$$
 P^{\otimes 2}  \left[ f(N_n) g\big( n^{-1/2}S_n
\big) \tilde
  g\big( n^{-1/2}\tilde S_n
\big)
\right]
\to
P^{\otimes 2}   [ f(N_\8)] \times P[ g(Z_1)] \times  P[ \tilde g(\widetilde Z_1)].
$$
Since $N_\8$  is geometrically distributed with parameter $1-\pi_d$ (see \eqref{eq-distri-Ninfty}), this ends the proof of
the lemma.
\qed
\begin{remark}
We could have taken another route to prove the lemma. The couple $n^{-1/2}(S_n, \tS_n)$ converges (mixing) to the Gaussian vector $(Z_1, \tilde Z_1)$, see e.g. Theorem 2
in~\cite{AldousEagleson78}. On the other hand $N_n \to N_\8$ a.s.
From the  mixing consequence  that we mentioned above Theorem \ref{th:clt} (which remains valid for random variables with values in $\R^d$)   it follows
  the convergence in distribution
  of $\big(n^{-1/2}S_n, n^{-1/2} \tS_n, N_\8\big)$ to $(Z_1, \tilde Z_1, N)$.
It is not difficult to see that the  sequence  $(n^{-1/2}S_n, n^{-1/2} \tS_n, N_\n)$ has the same limit  in distribution.
However, we have given the above proof, which is short and instructive, for the sake of completeness.
\end{remark}

\noindent
{\bf
Proof of Lemma \ref{lem:moments}}. (a) The "only if" part is evident by Fatou's lemma since
 $$  \Big| \frac{ S_n -\tilde S_n}{\sqrt n}  \Big|^{-a}  {\bf 1}_{\{  S_n - \tilde S_n \neq 0 \}} \rightarrow |Z|^{-a} \qquad \mbox{ in distribution }
 $$
  and
$E |Z|^{-a} < \infty $ if and only if $a<d$.

Let's show the "if" part.
Since the $L^a$-norm is increasing in $a$, we only need to prove the finiteness of the $\limsup$ for $0<a<d $ sufficiently close to $d$. So we fix $a \in (d-1, d)$. By the local central limit theorem
(see e.g. \cite[Theorem 1.2.1, p.14]{Lawler})
 we know that
$$ p_n(x) : = P^{\otimes 2}  (S_n - \tilde S_n =x ) = P(S_{2n} = x) $$
satisfies
\begin{equation} \label{eqn:pn}
 |p_n(x) - \bar p_n (x) | \leq c_2 n^{-d/2} |x|^{-2},
 \mbox{ with }
\bar p_n(x) =  c_1 n^{-d/2} \exp \{ - \frac{d |x|^2}{4n}\},
\end{equation}
where  $c_1, c_2>0$ are constants.
(In fact we have $c_1= \CstF$, with $\CstF$  defined by \eqref{def:CstF}.)
Denote by $I_n$  the integral in (\ref{eqn:moments}). Then
\begin{equation} \label{eqn:In}
 I_n \leq \sum_{x\neq 0} \big(\frac{|x|}{\sqrt n}\big)^{-a} \bar p_n(x)
 + \sum_{x\neq 0} \big(\frac{|x|}{\sqrt n}\big)^{-a} c_2 n^{-d/2} |x|^{-2} := I_{n,1} + I_{n,2}.
\end{equation}
In the following to avoid sums over non-integer valued numbers, we shall use the integer valued $L^1$ norm
$\| x \|_1 = |x_1| + \cdots + |x_d| $ for $x = (x_1, \cdots, x_d) \in \Z^d$, instead of the Euclidean  norm $|x|$, and the elementary  inequality  $ \|x\|_1/d \leq |x| \leq \| x \|_1$ valid for all $x\in \Z^d$.

For the first sum in (\ref{eqn:In}), we have, for some constant $c_3>0$,
\begin{eqnarray}
I_{n,1} & \leq &  c_1 d^a  n^{-d/2}
   \sum_{x\neq 0} \Big(\frac{\|x\|_1}{\sqrt n}\Big)^{-a} \exp \{ - \frac{\|x\|_1^2}{4dn}\}
       \nonumber \\
       & = &  c_1 d^a   n^{-d/2}
    \sum_{r \geq  1} \sum_{\|x\|_1=r} \Big(\frac{r}{\sqrt n}\Big)^{-a} \exp \{ - \frac{r^2}{4dn}\}
       \nonumber \\
 &\leq & c_1 \, c_3 \, d^a  n^{-d/2} \sum_{r \geq  1}  r^{d-1} \Big(\frac{r}{\sqrt n}\Big)^{-a}  \exp \{ - \frac{r^2}{4dn}\}
     \nonumber \\
     &=&   \frac{  c_1 \, c_3 \, d^a }{\sqrt n}   \sum_{r \geq  1} \Big(\frac{r}{\sqrt n}\Big)^{-(a-d+1)} \exp \{ - \frac{1}{4d} (\frac{r}{\sqrt n})^2\} , \nn
       \end{eqnarray}
       where the next to last step holds as the number of
      $ x = (x_1, ..., x_d)\in \zd$ with $|x_1| + \cdots + |x_d|=r$ is bounded by  $ 2(2r+1)^{d-1} \leq c_3 r^{d-1}$ (notice that each coordinate satisfies $|x_i| \leq r$, and while the first $d-1$ are chosen, the absolute value of the  last coordinate is determined by the equation, so that the last coordinate  has at most $2$ possibilities).  As $0< a-d+1 < 1 $, we have, for all $n\geq 1$,
      \begin{eqnarray}
 I_{n,1} & \leq &  \frac{c_1\, c_3  \, d^a  }{\sqrt n}  \Big( \sum_{1\leq r \leq \sqrt n}
        (\frac{r}{\sqrt n})^{-(a-d+1)}
        + \sum_{r > \sqrt n} \exp \{ -   \frac{1}{4d} (\frac{r}{\sqrt n})^2 \} \Big)
        \nonumber \\
       & \leq &
       c_1 \, c_3 \, d^a \Big( \sum_{1\leq r \leq \sqrt n}  \int_{\frac{r-1}{\sqrt n}}^\frac{r}{\sqrt n} x^{-(a-d+1)} dx
                +  \sum_{r > \sqrt n}  \int_{\frac{r-1}{\sqrt n}}^\frac{r}{\sqrt n}  \exp \{ -  \frac{x^2}{4d}\}   dx \Big)
                \nonumber \\
       & \leq & c_4:=  c_1 \, c_3 \, d^a  \Big(\int_0^1    x^{-(a-d+1)} dx
           + \int_0^\infty \exp \{ -\frac{x^2 }{4d} \} dx \Big)
             < \infty. \nn
                \end{eqnarray}
Similarly, for the second sum in (\ref{eqn:In}), we have,
\begin{eqnarray}
I_{n,2}
& \leq & c_2 \, c_3 \, d^{a+2} n^{-d/2} \sum_{r \geq 1}
     \Big(  \frac{r}{\sqrt n}\Big)^{-a}  r^{-2} r^{d-1} \nonumber \\
   &=& c_2 \, c_3 \, d^{a+2}   n^{(a-d)/2} \sum_{r\geq 1} r^{-(a+3-d)} \nn
   \end{eqnarray}
   which tends to $0$ as $n\rightarrow \infty$ (since $d\!-\!1 \!<\! \!a\!<\!d$).
This ends the proof of Part (a) of Lemma \ref{lem:moments}.

(b) Let $G_n=G_n(S_n - \tilde S_n) :=   n^{\frac{d-2}{2}} G(S_n - \tilde S_n)$ and $\lambda_2'=(1+\delta)\lambda_2$.
Recall from (\ref{eq:greenest}) that $G(x)\leq c_5 |x|^{2-d}$. Then

\begin{eqnarray}
 P^{\otimes 2}  \left[ (e^{\lambda_2 N_n } G_n )^{1+\delta} \right] &=& \sum_z
 P^{\otimes 2}    \left[ e^{\lambda_2' N_n } ; S_n-\tS_n=z \right]  G_n(z)^{1+\delta}
 \nn
 \\ &\leq& c_5 \sum_{z \neq 0}
 P^{\otimes 2}    \left[ e^{\lambda_2' N_n } ; S_n-\tS_n=z \right]  \Big|n^{-1/2}z\Big|^{-(d-2)(1+\delta)}  \nn
\\ &&\qquad  +\, c_5
 P^{\otimes 2}    \left[ e^{\lambda_2' N_n } ; S_n-\tS_n=0 \right]  n^{(d-2)(1+\delta)/2} . \label{eq:pac1}
\end{eqnarray}
Denoting by
 $L_n=\sup\{j =0,\ldots n: S_j-\tS_j=0\}$ and $T=\inf \{j \geq 1: S_j-\tS_j=0\}$
 the last (before $n$) and the first
hitting times of $0$, we will use that, under Condition {\bf (L2)},
\begin{equation}
\label{eq:GB1}
 P^{\otimes 2}    \left[ e^{\lambda_2' N_n } ; S_n-\tS_n=0 \right] \sim e^{\lambda_2'} [1-e^{\lambda_2'}\pi_d]^{-2} P^{\otimes 2}   [T=n] \qquad  {\rm as \ } n \to \8,
\end{equation}
see Theorem 2.2 (case 2) in \cite{GB}.
Since $P^{\otimes 2}   [T=n] \leq P^{\otimes 2}  (S_n-\tS_n=0) = O(n^{-d/2})$,   this implies that when $(d-2)(1+\delta)<d$,
                      the last term in
 (\ref{eq:pac1}) vanishes as $n\to \infty$, and yields, for $z \neq 0$,
\begin{eqnarray}
P^{\otimes 2}    \left[ e^{\lambda_2' N_n } ; S_n-\tS_n=z \right]  &=&
\sum_{j=0}^{n} P^{\otimes 2}    \left[ e^{\lambda_2' N_n } ; S_n-\tS_n=z , L_n=j \right]
\nn \\ &\stackrel{\rm Markov}{=}&
\sum_{j=0}^{n} P^{\otimes 2}    \left[ e^{\lambda_2' N_{j} } ; S_{j}\!-\!\tS_{j}=0 \right]
P^{\otimes 2}    \left[  S_{n-j}\!-\!\tS_{n-j}=z , L_{n-j}=0 \right]
\nn  \\ &\stackrel{(\ref{eq:GB1})}{\leq}&
c_6  \sum_{j=0}^{n} P^{\otimes 2}   [T=j] P^{\otimes 2}    \left[  S_{n-j}\!-\!\tS_{n-j}=z , L_{n-j}=0 \right]
\nn  \\ &{\leq}&
c_6 \sum_{j=0}^{n} P^{\otimes 2}   [T=j] P^{\otimes 2}    \left[  S_{n-j}\!-\!\tS_{n-j}=z \right] \nn \\
&\leq & c_6 P^{\otimes 2}    \left[  S_{n}\!-\!\tS_{n}=z \right]. \nn
\end{eqnarray}
Inserting this in (\ref{eq:pac1}), we obtain for large $n$,
\begin{eqnarray}
 P^{\otimes 2}  \left[ (e^{\lambda_2 N_n } G_n )^{1+\delta} \right] &=& c_5 c_6 \left[ 1+ \sum_{z \neq 0}
 P^{\otimes 2}    \left[ S_n-\tS_n=z \right]  \Big|n^{-1/2}z\Big|^{-(d-2)(1+\delta)} \right], \nn
\end{eqnarray}
which, according to (\ref{eqn:moments}), is bounded
provided that $(d-2)(1+\delta)<d$.
 \qed

\section{Proof of the Central Limit Theorem} \label{sec:3}

In this section we give the proof of Theorem \ref{th:clt}.
Our proof is based on the following central limit theorem for infinite martingale arrays, which is a slight extension of Corollaries 3.1 and 3.2 of the book by  Hall and Heyde \cite{Hall-Heyde} (pp. 58-59 and p.64), but we could not find it in the literature.

\begin{lemma} \label{lem:martingale-clt}
For $n\geq 1$, let $ \{ (S_{n,i}, {\cal F}_{n,i}): i\geq 0 \}$ be a martingale defined on a probability space
$(\Omega, {\cal F}, P)$, with $S_{n,0}=0$ and
\begin{equation} \label{cond-L2bdd}
    \sup_{n,i \geq 1} E S_{n,i}^2 < \infty.
\end{equation}
Let $X_{n,i} = S_{n,i} - S_{n,i-1}, i\geq 1$ be the martingale differences, and  $  S_{n,\infty} = \lim_{i\rightarrow \infty} S_{n,i}$ be the a.s. limit of $(S_{n,i}, i\geq 0)$. Suppose that:
\begin{enumerate}
  \item[(a)] the conditional variance converges in probability:  for a real random variable $V \in [ 0, \infty)$,
     \begin{equation} \label{cond-cond-variances}
        V^2_{n,\infty} : = \sum_{i=1}^\infty E (X_{n,i}^2 | {\cal F}_{n, i-1}) \longrightarrow V^2   \mbox{ in probability; }
        \end{equation}
  \item[(b)] the conditional Lindeberg condition holds :
   \begin{equation} \label{cond-cond-lind}
  \forall \varepsilon >0, \quad
    \sum_{i=1}^\infty E (X_{n,i}^2  {\bf 1}_{ |X_{n,i} | > \varepsilon} | {\cal F}_{n, i-1} ) \longrightarrow 0
      \mbox{ in probability; }
  \end{equation}
  \item[(c)] the $\sigma-fields$ are nested: $ {\cal F}_{n, i} \subset {\cal F}_{n+1, i}$ for all $n,i\geq 1$.
\end{enumerate}
Then
\begin{equation}\label{conclu-mart-clt}
  S_{n,\infty} \longrightarrow V \, G  \quad \mbox{ in distribution }
\end{equation}
where $G$ is a Gaussian variable with law $N(0,1)$ and independent of $V$; if additionally  $V \neq 0$ a.s., then
\begin{equation}\label{conclu-mart-clt2}
  \frac{ S_{n,\infty}}{V_{n,\infty}} \longrightarrow G  \quad \mbox{ in distribution. }
\end{equation}
Moreover, the convergence in \eqref{conclu-mart-clt} is stable, and the convergence in \eqref{conclu-mart-clt2} is mixing.
\end{lemma}

 Lemma \ref{lem:martingale-clt} reduces to Corollaries 3.1 and 3.2 in \cite{Hall-Heyde}   for a triangular array of martingales differences, that is,  when $X_{n,i} = 0$ for all $i >k_n$,
  for some sequence of integers $k_n$ increasing to $\infty$.
 As in the case of a triangular array,
 if $V$ is measurable with respect to each ${\cal F}_{n,i}$ for $n,i \geq 1$ (e.g. when $V$ is a constant), then the nested condition (c) can be removed, but the convergence \eqref{conclu-mart-clt} may no longer be stable, and the  convergence
\eqref{conclu-mart-clt2} may  no longer be mixing (see the remarks in p.59 and p. 64 of \cite{Hall-Heyde} for a triangular array).
 Lemma \ref{lem:martingale-clt} can be extended in a clear way to two-sided martingale arrays $ \{ (S_{n,i}, {\cal F}_{n,i}): -\infty <i < \infty \}$; for a  version using conditions and norming in terms of $ \sum_i X_{n,i}^2 $, see Theorem 3.6 of \cite{Hall-Heyde} (p.77).

\medskip
\noindent
{\bf Proof of Lemma \ref{lem:martingale-clt}}.  We will see that Lemma \ref{lem:martingale-clt} can be obtained from the corresponding result for a triangular array  of martingales. Let $k_n$ be positive integers increasing to $\infty$ such that
$$ E (S_{n,\infty} - S_{n,k_n} )^2 = E \sum_{i>k_n} X^2_{n,i} \rightarrow 0. $$
Then
$$ V_{n,k_n}^2 :=  \sum_{i=1}^{k_n} E (X_{n,i}^2 | {\cal F}_{n, i-1}) \rightarrow V^2   \mbox{ in probability} $$
since $ V_{n,\infty}^2 - V_{n,k_n}^2 =  \sum_{i>k_n} E (X_{n,i}^2 | {\cal F}_{n, i-1}) \rightarrow 0 $ in $L^1$. Clearly, by condition (b), the conditional Lindeberg condition for the triangular array $\{(X_{n,i}, {\cal F}_{n,i}): 0\leq i\leq k_n \}$ holds, that is, \eqref{cond-cond-lind} holds with $\sum_{i=1}^\infty$ replaced by $\sum_{i=1}^{k_n}$. Therefore, by
Corollaries 3.1 and 3.2 of  \cite{Hall-Heyde} (pp.58-59 and p. 64),
\begin{equation}\label{conclu-mart-cltb}
  S_{n,k_n} \rightarrow V \, G  \quad \mbox{ in distribution }
\end{equation}
and
\begin{equation}\label{conclu-mart-clt2b}
  \frac{ S_{n,k_n}}{V_{n,k_n}} \rightarrow G  \quad \mbox{ in distribution. }
\end{equation}
Since $ S_{n,\infty} - S_{n,k_n} \rightarrow 0$ in $L^2$ and hence in probability, \eqref{conclu-mart-cltb}
implies \eqref{conclu-mart-clt}. As $ V^2_{n,\infty} - V^2_{n,k_n} \rightarrow 0$  in probability (in fact in $L^1$), when $V>0$ a.s. we have $ V^2_{n,\infty} / V^2_{n,k_n} \rightarrow 1$  in probability.
Therefore \eqref{conclu-mart-clt2b}
implies \eqref{conclu-mart-clt2}. \qed

Lemma \ref{lem:martingale-clt} is well suited for  studying the rate of convergence of  a martingale, as shown in the following

\begin{corollary} \label{coro-martingale-conv-rate}
Let $ \{ (S_{i}, {\cal F}_{i}): i\geq 0 \}$ be a martingale defined on a probability space
$(\Omega, {\cal F}, P)$, with $S_{0}=0$ and
$
    \sup_{i \geq 1} E S_{i}^2 < \infty.
$
Let $X_{i} = S_{i} - S_{i-1}, i\geq 1$ be the martingale differences,   $  S_{\infty} = \lim_{i\rightarrow \infty} S_{i}$ be the a.s. limit of $(S_{i})$, and let
$$v^2_n = E(S_\infty - S_n)^2  = E \sum_{i=n+1}^\infty X_i^2.$$
Suppose that $v_n>0$ and that:
\begin{enumerate}
  \item[(a)] the conditional variance converges in probability:  for a real random variable $V \in [ 0, \infty)$,
     \begin{equation} \label{cond-cond-variances-coro}
        V^2_{n} : = \frac{1}{v_n^2} \sum_{i=n+1}^\infty E (X_{i}^2 | {\cal F}_{i-1}) \rightarrow V^2   \mbox{ in probability; }
        \end{equation}
  \item[(b)] the conditional Lindeberg condition holds :
   \begin{equation} \label{cond-cond-lind-coro}
  \forall \varepsilon >0, \quad
    \frac{1}{v_n^2} \sum_{i=n+1}^\infty E (X_{i}^2  {\bf 1}_{ |X_{i} | > \varepsilon v_n } | {\cal F}_{i-1} ) \rightarrow 0
      \mbox{ in probability. }
      \end{equation}
\end{enumerate}
Then
\begin{equation}\label{conclu-mart-clt-coro}
  \frac{S_{\infty}- S_n}{v_n}  \rightarrow V \, G  \quad \mbox{ in distribution, }
\end{equation}
where $G$ is a Gaussian variable with law $N(0,1)$ and independent of $V$; if additionally  $V  \neq 0$ a.s., then
\begin{equation}\label{conclu-mart-clt2-coro}
  \frac{ S_{\infty} - S_n}{V_{n}} \rightarrow G  \quad \mbox{ in distribution. }
\end{equation}
Moreover, the convergence in \eqref{conclu-mart-clt-coro} is stable, and the convergence in \eqref{conclu-mart-clt2-coro} is mixing.
\end{corollary}

\medskip
\noindent
{\bf Proof.} Corollary \ref{coro-martingale-conv-rate} is an immediate consequence of Lemma \ref{lem:martingale-clt}
 applied to  $X_{n,i} = X_{n+i}$  for $i \geq 1$,  $ {\cal F}_{n,i} = {\cal F}_{n+i} $ for $i \geq 0$, and $X_{n,0}=0$.  \qed

\medskip

\noindent
{\bf Proof of Theorem \ref{th:clt}}. We can use Corollary \ref{coro-martingale-conv-rate} with
$ 1/ \| W-W_n \|_2  $ for the norming, but, since $\| W-W_n \|_2  \sim \sigma n^{-(d-2)/4}$, we prefer to use the
more explicit norming $ n^{(d-2)/4}$, together with the spirit of the proof of Corollary \ref{coro-martingale-conv-rate}.
So we rely on the decomposition
\begin{equation} \label{eq:mgdec}
  n^{ \frac{d-2}{4}} (W-W_n) =  n^{ \frac{d-2}{4}} \sum_{k=n}^\infty  D_{k+1},
\end{equation}
where 
\begin{equation}\nn
D_{k+1}= W_{k+1}-W_k , \quad k \geq n,
\end{equation}
forms a sequence of martingale differences.
To prove Theorem \ref{th:clt}, by \eqref{eq:mgdec} and
 Lemma \ref{lem:martingale-clt}
 applied to  $X_{n,i} = n^{ \frac{d-2}{4}} D_{n+i}$  for $i \geq 1$,  $ {\cal F}_{n,i} = \cG_{n+i} $ for $i \geq 0$ (recall \eqref{def-algGn}), and $X_{n,0}=0$,
 it suffices  to prove that :
\begin{enumerate}
\item[(a)]  the following convergence about the conditional variance holds  :
\begin{equation} \label{cond:var}
 s_n^2:= n^{\frac{d-2}{2}} \sum_{k\geq n}  \E_k D_{k+1}^2   \rightarrow  \sigma_1^2 W^2  \quad \mbox{ in probability, }
\end{equation}
where $\E_k (\cdot) = \E ( \cdot | \cG_k)$ denotes the conditional expectation given $\cG_k$;
\item[(b)]    the following  Lindeberg condition holds :
\begin{equation}  \label{cond:lind}
 \forall \eps>0, \quad n^{ \frac{d-2}{2}} \sum_{k\geq n} \,  \E_k  \left(  D_{k+1}^2
    {\bf 1}_{  \{ n^{ \frac{d-2}{4}}   |D_{k+1}|  >  \varepsilon   \}  } \right) \to 0 \quad \mbox{ in probability.}
\end{equation}
\end{enumerate}
Actually, by Lemma \ref{lem:martingale-clt}, from \eqref{cond:var} and \eqref{cond:lind} we conclude that
\eqref{clt2bis1} and  \eqref{clt2bis} hold with the norming  $1/W_n$ in \eqref{clt2bis} replaced by $1/W$.
As $W_n/W \to 1$ a.s. (and thus in probability),  we can change the factor    $1/W$  to $1/W_n$ without changing the convergence in distribution.

To show the convergence \eqref{cond:var} of the conditional variance,
we will prove
 in the next section
the following
\begin{lemma} \label{lem:instr}
There exists $\beta_0>0$ such that for $|\beta|<\beta_0$ and $\sigma_1$ from \eqref{eq:sigma1}, we have, as $n \to \8$,
\begin{eqnarray}
 \label{lem:conv-L4}
 \E ( W_n  -  W)^4  &\longrightarrow& 0,   \\
\label{lem:instr4}
\E s_n^4 -  \sigma_1^4 \E W_n^4 &\longrightarrow& 0,\\
\E (s_n^2 W_n^2) -  \sigma_1^2 \E W_n^4 &\longrightarrow& 0.  \label{lem:instr3}
\end{eqnarray}
\end{lemma}

It follows  from Lemma \ref{lem:instr} that
\begin{eqnarray}\nn
\E (s_n^2 -  \sigma_1^2 \E W_n^2)^2 &=& \E s_n^4  -2  \sigma_1^2  \E (s_n^2 W_n^2) + \sigma_1^4 \E W_n^4 \\
& \longrightarrow& 0 . \nn
\end{eqnarray}
Therefore, $s_n^2 -  \sigma_1^2  W_n^2 \to 0$ in $L^2$.
 As $W_n^2 \to W^2$ in $L^2$, it follows that $s_n^2 \to   \sigma_1^2  W^2 $ in $L^2$.  We thus obtain \eqref{cond:var}.
\medskip

To show  Lindeberg's condition \eqref{cond:lind}, we will prove the following convergence rate
of $\E D^4_{k+1}$.

\begin{lemma} \label{lem:rate-conv-Dk} For any $q>1$, when $|\beta|>0$ is small enough, we have
\begin{equation} \label{eq-rate-conv-Dk}
\E D^4_{k+1} = O (k^{-d/q}), \qquad k\geq 1.
\end{equation}
\end{lemma}
The proof of the Lemma is postponed to the end of the paper.

Note that a sufficient condition for  Lindeberg's condition \eqref{cond:lind} to hold is clearly
\begin{equation}\nn 
n^{  (d-2)  } \sum_{k\geq n}   \E_k D_{k+1}^{4}
 \rightarrow 0 \quad \mbox{ in probability,  }
\end{equation}
which is implied by
\begin{equation} \label{cond-lind2bis}
n^{  (d-2)  } \sum_{k\geq n}   \E D_{k+1}^{4}
 \rightarrow 0.
\end{equation}
By Lemma \ref{lem:rate-conv-Dk}, when $|\beta|>0$ is small enough, the left-hand side of \eqref{cond-lind2bis} is smaller than
 $ c \, n^{d-1 - d/q  }$ for some constant $c>0$, which tends to $0 $ by taking  $1< q < d/(d-1)$,
so that Lindeberg's condition \eqref{cond:lind} holds. (In fact one can check that it suffices to take $|\beta| < \beta_0$ with $\beta_0 >0$ determined in Lemma \ref{lem:vargas}.)
 This ends the
proof of Theorem \ref{th:clt}, using  Lemmas \ref{lem:instr} and \ref{lem:rate-conv-Dk} whose proofs will be given in the next section.  \qed

\begin{remark}  The convergence \eqref{clt2bis} can also be proved using the decomposition
\begin{equation} \nn \label{eq:mgdec-bis}
   \frac{n^{ \frac{d-2}{4}} (W-W_n)}{W_n}  =   \frac{n^{ \frac{d-2}{4}} }{W_n}  \sum_{k=n}^\infty  D_{k+1},
\end{equation}
where  $(n^{\frac{d-2}{4}}D_{k+1}/W_n, k\geq n)$ remains a sequence of martingale differences.
But proving \eqref{clt2bis1} requires a different route.
\end{remark}

\medskip

We end this section with a {\bf warm-up calculation}: we recover the value of  $\E W_n^2$ and (\ref{eq:dist2}) from the martingale decomposition \eqref{eq:mgdec},
i.e., using the conditional variance.
This calculation is instructive and it will be useful in the forthcoming computations. Write for short
 $$ h_k(S)  = \sum_{i=0}^{k-1}  [ \beta \eta(i, S_i) -  \lambda (\beta)].
$$
Since  $W_{k+1} - W_k = P e^{h_k(S)} (e^{ \beta \eta (k,S_k) - \lambda (\beta)} -1)$,
$$ (W_{k+1} - W_k)^2  = P^{\otimes 2} \left[ e^{h_k(S)} (e^{ \beta \eta (k,S_k) - \lambda (\beta)} -1)
                                      e^{h_k(\tilde S )} (e^{ \beta \eta (k,\tilde S_k) - \lambda (\beta)} -1) \right],
                                      $$
using Fubini's theorem
we have
\begin{equation}
 \E_k D_{k+1}^2 = \Cstd P^{\otimes 2} e^{h_k(S)} e^{h_k(\tilde S )}
           {\bf 1}_{ \{ S_{k} = \tilde S_{k} \}}
 \end{equation}
    with
    \begin{equation} \label{def:kappa2}
\Cstd= \Cstd(\beta)= e^{\lambda_2}-1 .
\end{equation}
We compute $s_n^2$ from its definition in \eqref{cond:var}:
\begin{equation}
\label{def:sigman}
s_n^2 = \Cstd n^{ \frac{d-2}{2}} \sum_{k\geq n}  P^{\otimes 2} e^{h_k(S)} e^{h_k(\tilde S )}
           {\bf 1}_{S_{{k}} = \tilde S_{{k}}}.
           \end{equation}
Observe that
 \begin{equation}
\label{eq:bki52}
   {\bf 1}_{\{S_{{k}} = \tilde S_{{k}}\}} = \frac{e^{\lambda_2 (N_{k+1}-N_k)}-1}{\Cstd}.
\end{equation}

We now check that, with  $s_n^2$ from \eqref{def:sigman} and  $\s^2$ be defined by (\ref{eq:sigma2}),
$$ \E s_n^2 \rightarrow \s^2.   $$
This is a remake of \eqref{eq:dist2cv}, but, as we will see,  from a different route.
By Fubini we have
\begin{eqnarray} \label{eq-esp-sn}
 \E s_n^2 & = &\Cstd n^{ \frac{d-2}{2}}  \sum_{k\geq n}  P^{\otimes 2} \E e^{h_k(S)} e^{h_k(\tilde S )}
                             {\bf 1}_{S_{{k}} = \tilde S_{{k}}} \nonumber \\
          & =&  n^{ \frac{d-2}{2}}    P^{\otimes 2}  \sum_{k\geq n} e^{\lambda_2 N_k}
             \Cstd {\bf 1}_{S_{{k}} = \tilde S_{{k}}}  \nonumber \\
          & =&  n^{ \frac{d-2}{2}}    P^{\otimes 2} \sum_{k\geq n} e^{\lambda_2 N_k}
               (e^{\lambda_2 (N_{k+1}-N_{k})} -1)  \qquad ({\rm by\ } \eqref{eq:bki52}) \nonumber \\
          & =&  n^{ \frac{d-2}{2}}    P^{\otimes 2} ( e^{\lambda_2 N_\infty}
                - e^{\lambda_2 N_n} )    \qquad \qquad \qquad ({\rm telescopic \ sum})\nonumber \\
          & =&  n^{ \frac{d-2}{2}}    P^{\otimes 2} [P^{\otimes 2} \left( e^{\lambda_2 (N_\infty - N_n)} - 1\right) | \cF_n )
                e^{\lambda_2 N_n}  ]    \nonumber \\ \label{eq:j08}
           & =&  n^{ \frac{d-2}{2}}   P^{\otimes 2} [  F (S_n-\tilde S_n )
                e^{\lambda_2 N_n}  ],    \nn
\end{eqnarray}
where we can express
\begin{eqnarray}
F(x)
 &=&  P_{0,x}^{\otimes 2} (e^{\lambda_2 N_\infty} -1) \nonumber \\
  &=&  P_{0,x}^{\otimes 2} ((e^{\lambda_2 N_\infty} -1)  {\bf 1}_{\tau < \infty}) \nonumber \\
 &=&  P_{0,x}^{\otimes 2} (\tau < \infty)  P_{0,0}^{\otimes 2} ( e^{\lambda_2 N_\infty}  - 1) \nonumber \\
 &=& \frac{G(x)} { G(0)} (\E (W^2) - 1).  \label{def-F}
\end{eqnarray}
Therefore, by the same argument as in \eqref{eq:grosqnu},
 we derive from the last 2 formulas, 
\begin{eqnarray} \label{eq-esp-sn2}
 \E s_n^2 & = &  \frac{\Var(W) }{G(0)}   P^{\otimes 2} [   n^{ \frac{d-2}{2}}  G (S_n-\tilde S_n )
                e^{\lambda_2 N_n}  ] \nonumber \\
                & \rightarrow & \frac{\Cst }{G(0)}  \Var(W)  E \frac{1}{|Z| ^{d-2} } \times P^{\otimes 2} e^{\lambda_2 N_\8} \nonumber \\
                &=&{\Cst   {\rm Var} (W)}{(1-\pi_d)} \Zd  \times \E W^2 \nn
                \end{eqnarray}
                which is equal to $\sigma^2$.


%
%

\section{Proof of the instrumental   lemmas} \label{sec:4}


In this section we give the proofs of Lemmas \ref{lem:instr} and \ref{lem:rate-conv-Dk}.


\medskip
\noindent
{\bf Proof of Lemma \ref{lem:instr}.}  We first give the {\bf proof of \eqref{lem:conv-L4}}.
For 4 independent paths $S^{(i)}$,
$i =1,\ldots 4$, we need to count the number of intersections of 2, 3 or four of them. Denote by $E_4$ the set of elements of $(\Z^d)^4$
with all 4 elements equal, by $E_3$ the set of those with 3 equal lattice sites and a different fourth, by $E_{2,2}$ the set of those with 2 pairs of equal lattice sites but the two are different, and by $E_{2,0}$ the set of those with one pair being equal and different from the two other ones. Let $A$  be the corresponding set of indices,
$A=\{4, 3, (2,2), (2,0)\}$, and define
\begin{equation}
\nn
\lambda_m=\lambda(m\beta)-m \lambda(\beta),\quad  m \geq 2, \qquad
N_{n,k}^{(a)} = \sum_{t=n}^{k-1}
{\bf 1}_{ \{(S^{(1)}_t,S^{(2)}_t,S^{(3)}_t,S^{(4)}_t) \in E_a\}}, \quad a \in A,
\end{equation}
\begin{equation}
\nn
\Sigma_{n,k}=\lambda_4 N_{n,k}^{(4)}+ \lambda_3 N_{n,k}^{(3)}+ 2\lambda_2 N_{n,k}^{(2,2)}+
\lambda_2 N_{n,k}^{(2,0)} , \qquad 0 \leq n \leq k.
\end{equation}
                 Then it is elementary to check that
\begin{equation} \label{eq:alg4}
\E \exp \{ \sum_{i=1}^4 h_n(S^{(i)})\} = \exp \Sigma_{0,n},    \qquad \E W_n^4 = P^{\otimes 4}\left(e^{\Sigma_{0,n} }\right) .
\end{equation}
For all pair $S, \tS$ of paths and all $0 \leq n \leq k$, put $N_{n,k}(S, \tS)= \sum_{i=n}^{k-1} {\bf 1}_{S_i=\tS_i}$. Note that, for all $a \in A$,
\begin{equation} \label{eq:alg4bds}
N_{n,k}^{(a)} \leq \sum_{ 1\leq i < j \leq 4} N_{n,k}(S^{(i)},S^{(j)}),
\end{equation}
and that $\lambda_m \searrow 0$ as $|\beta| \searrow 0$. Thus, taking $|\beta|>0$ small, we have uniform integrability of $(W_n^4)_n$, and then
\begin{equation}\label{eq:int1}
\sup_n \E W_n^{4}= \E W^4 < \8 .
\end{equation}
This gives \eqref{lem:conv-L4}.

We next give the {\bf proof of \eqref{lem:instr4}}.
We compute
\begin{eqnarray}
\E s_n^4 &=&
\nn
\Cstd^2 n^{ d-2}  \sum_{k, \ell \geq n}  P^{\otimes 4} \left[ \E e^{h_k(S^{(1)})+h_k(S^{(2)})+ h_\ell(S^{(3)})+h_\ell(S^{(4)})}
                             {\bf 1}_{S^{(1)}_{{k}} =  S^{(2)}_{{k}}}  {\bf 1}_{S^{(3)}_{{\ell}} =  S^{(4)}_{{\ell}}} \right]\nonumber \\  &=&  T^<(n) + T^=(n) + T^>(n) , \nn
\end{eqnarray}
where $T^=(n)$ [resp. $T^<(n)$, resp. $T^>(n)$]  is the contribution in the sum of the terms with $k=\ell$, [resp. $k<\ell$, resp.
$k>\ell$]. 
By symmetry, $T^<(n)=T^>(n)$. We calculate $T^<(n) $:
\begin{eqnarray} \nn
T^<(n)       & =& \Cstd^2   n^{{d-2}}    P^{\otimes 4}  \left[\sum_{\ell> k\geq n} e^{\Sigma_{0,k}+\lambda_2 N_{k,\ell}(S^{(3)},S^{(4)})}  {\bf 1}_{S^{(1)}_{{k}} = S^{(2)}_{{k}}} {\bf 1}_{S^{(3)}_{{\ell}} = S^{(4)}_{{\ell}}}\right]  \nonumber \\
   & =& \Cstd   n^{{d-2}}    P^{\otimes 4} \left[ \sum_{ k\geq n} e^{\Sigma_{0,k}} {\bf 1}_{S^{(1)}_{{k}} = S^{(2)}_{{k}}} \sum_{\ell>k} e^{ \lambda_2 N_{k,\ell}(S^{(3)},S^{(4)})} \Cstd  {\bf 1}_{S^{(3)}_{{\ell}} = S^{(4)}_{{\ell}}} \right]  \nonumber \\
         &\stackrel{\eqref{eq:bki52}}{ =}&  \Cstd n^{ {d-2}}    P^{\otimes 4}  \left[ \sum_{k\geq n}
 e^{\Sigma_{0,k}} {\bf 1}_{S^{(1)}_{{k}} = S^{(2)}_{{k}}} \sum_{\ell>k} e^{ \lambda_2 N_{k,\ell}(S^{(3)},S^{(4)})}
\left(  e^{ \lambda_2 N_{\ell,\ell+1}(S^{(3)},S^{(4)})}-1\right) \right]
  \nonumber \\
     &=&  \Cstd n^{ {d-2}}    P^{\otimes 4}  \left[ \sum_{k\geq n}
 e^{\Sigma_{0,k}} {\bf 1}_{S^{(1)}_{{k}} = S^{(2)}_{{k}}} \sum_{\ell>k}\left( e^{ \lambda_2 N_{k,\ell+1}(S^{(3)},S^{(4)})}-   e^{ \lambda_2 N_{k,\ell}(S^{(3)},S^{(4)})}\right)
\right]  \nonumber \\
     &=&  \Cstd n^{ {d-2}}    P^{\otimes 4} \left[ \sum_{k\geq n}
 e^{\Sigma_{0,k}} {\bf 1}_{S^{(1)}_{{k}} = S^{(2)}_{{k}}}\left( e^{ \lambda_2 N_{k,\8}(S^{(3)},S^{(4)})}-   e^{ \lambda_2 N_{k,k+1}(S^{(3)},S^{(4)})}\right)
 \right] \nonumber \\
    &=&  \Cstd n^{ {d-2}}    P^{\otimes 4} \left[ \sum_{k\geq n}
 e^{\Sigma_{0,k}} {\bf 1}_{S^{(1)}_{{k}} = S^{(2)}_{{k}}}\left( e^{ \lambda_2 N_{k,\8}(S^{(3)},S^{(4)})}-  1
\right) \right] - T^=(n)
  \nonumber \\
     &\stackrel{\rm }{=}&  \Cstd n^{ {d-2}}    P^{\otimes 4}  \left[ \sum_{k\geq n}
 e^{\Sigma_{0,k}} {\bf 1}_{S^{(1)}_{{k}} = S^{(2)}_{{k}}}
  P^{\otimes 2}\left( e^{ \lambda_2 N_{k,\8}(S^{(3)},S^{(4)})}-  1 \big \vert \cF_k  \right) \right]
- T^=(n)
  \nonumber \\
     &\stackrel{\rm \eqref{def-F}}{=}&
\frac{{\rm Var}(W)}{G(0)}  \Cstd  \times  n^{ {d-2}}  P^{\otimes 4} \left[ \sum_{k\geq n}
 e^{\Sigma_{0,k}} {\bf 1}_{S^{(1)}_{{k}} = S^{(2)}_{{k}}}
 G(S^{(3)}_k-S^{(4)}_k)\right]
- T^=(n)
  \label{eq:clashantoine} \nn
  \\
  \nn
  &=&
   T^<_1(n) - T^= (n)\;,
\end{eqnarray}
which serves also as the definition of  $ T^<_1(n)$. For $0\leq m \leq n$ define
\begin{equation}
\nn
T^<_2(m,n) =   \frac{{\rm Var}(W)}{G(0)}  \Cstd
\times  n^{ {d-2}}  \sum_{k\geq n} P^{\otimes 4} \left[
 e^{\Sigma_{0,m}+ \lambda_2 N_{k-m, k}(S^{(1)}, S^{(2)})} {\bf 1}_{S^{(1)}_{{k}} = S^{(2)}_{{k}}}
 G(S^{(3)}_k-S^{(4)}_k)\right].
\end{equation}
It is sufficient to study the limit of $T^<_2(m,n)$ for large $n, m$, since we will prove the following
\begin{lemma} \label{lem:vargas}
Let $\beta_0>0$ be such that, for some $\eps>0$ and  all $|\beta| < \beta_0$,
 $\sup_n P^{\otimes 4} [  e^{(d+\eps)\Sigma_{0,n}} ] < \8$.
  Then,
$$
\lim_{n \to \8} \E  |T^=(n)| = \lim_{m \to \8} \limsup_{n \to \8} \E \big\vert T^<_1(n) - T^<_2(m,n)  \big\vert = 0 .
$$
\end{lemma}
\begin{remark}
The occurrence of the term $N_{k-m,k}$ may be surprising. It is reminiscent of a similar phenomenon in the local limit theorem for polymers
\cite{Si95, Va04}. The constraint $S_k^{(1)}=S_k^{(2)}$  makes likely intersections between
 $S^{(1)}$ and $S^{(2)}$
     just before time $k$, whereas it is likely that $S^{(3)}$ and $S^{(4)}$ are far from them and far apart one from another.
\end{remark}
With Lemma \ref{lem:vargas} we continue our proof. To analyse $T^<_2(m,n)$ we condition on the vectors
$$\cS^{(1,2)}_{m,k}=\cS^{(1,2)}=(S^{(i)}_t; i=1,2, t=0,\ldots m-1, {\rm \ and\ } t=k),$$
and $\cS^{(3,4)}=
(S^{(i)}_t; i=3,4, t=0,\ldots m-1)$, and use the independence of the paths $S^{(i)}$:
\begin{eqnarray}
T^<_2(m,n) \!\!\!\!
\nn &=
\frac{{\rm Var}(W)}{G(0)}  \Cstd n^{ {d-2}}  \sum_{k\geq n} P^{\otimes 4} \left[
 e^{\Sigma_{0,m}} {\bf 1}_{S^{(1)}_{{k}} = S^{(2)}_{{k}}}
 P^{\otimes 2} \big( e^{\lambda_2 N_{k-m, k}(S^{(1)}, S^{(2)})} \big \vert  \cS^{(1,2)}\big) \right. \\
 & \qquad \qquad  \qquad \qquad \qquad \qquad  \times \left.
 P^{\otimes 2} \big(  G(S^{(3)}_k-S^{(4)}_k) \big \vert \cS^{(3,4)} \big)
 \right] \nn \\
 \sim  &
 \frac{{\rm Var}(W) \Cst \Zd}{G(0)}  \Cstd n^{ {d-2}}  \sum_{k\geq n} \frac{1}{k^{\frac{d-2}{2}}} P^{\otimes 4} \left[
 e^{\Sigma_{0,m}} {\bf 1}_{S^{(1)}_{{k}} = S^{(2)}_{{k}}}
 P^{\otimes 2} \big( e^{\lambda_2 N_{k-m, k}(S^{(1)}, S^{(2)})} \big \vert  \cS^{(1,2)}\big) \right]  \label{eq:enfin}
\end{eqnarray}
as $n \to \8$, since for fixed $m$,  $k^{\frac{d-2}{2}} P^{\otimes 2} \big(  G(S^{(3)}_k-S^{(4)}_k) \big \vert \cS^{(3,4)} \big) \to \Cst \Zd$ uniformly as $k \to \8$ by  \eqref{eq:NYC2}, \eqref{eqn:moments}, \eqref{eqn:momentsG}.
Now,  on the event $\{S^{(1)}_{{k}} = S^{(2)}_{{k}}\} \bigcap \{S_{m-1}^{(i)}=s_{m-1}^{(i)}, i=1,2\}$, we have by time-reversal and Markov property,
\begin{eqnarray}
\nn
P^{\otimes 2} \big( e^{\lambda_2 N_{k-m, k}(S^{(1)}, S^{(2)})} \big \vert  \cS^{(1,2)}\big)
&=&
e^{-\lambda_2} P^{\otimes 2} \big( e^{\lambda_2 N_{0, m}} \big \vert  S^{(i)}_{k-m}=s^{(i)}_{m-1}, i=1,2\big)\\
\label{eq:levrard}
&\to& e^{-\lambda_2}  P^{\otimes 2} \big( e^{\lambda_2 N_{0, m}} \big) ,
\end{eqnarray}
as $k \to \8$. Thus, for fixed $m$ and $k \to \8$, it holds
\begin{eqnarray}
k^{\frac{d}{2}}
P^{\otimes 4} \left[
 e^{\Sigma_{0,m}} {\bf 1}_{S^{(1)}_{{k}} = S^{(2)}_{{k}}}
 P^{\otimes 2} \big( e^{\lambda_2 N_{k-m, k}} \big \vert  \cS^{(1,2)}\big) \right]
& \nn
\sim&
 e^{-\lambda_2} k^{\frac{d}{2}}
P^{\otimes 4} \left[
 e^{\Sigma_{0,m}} {\bf 1}_{S^{(1)}_{{k}} = S^{(2)}_{{k}}} \right] P^{\otimes 2} \big( e^{\lambda_2 N_{0, m}} \big) \\ \nn
 &\sim&
 e^{-\lambda_2} k^{\frac{d}{2}}
P^{\otimes 4} \left[
 e^{\Sigma_{0,m}}\right]  P^{\otimes 2} \big( {\bf 1}_{S^{(1)}_{{k}} = S^{(2)}_{{k}}} \big) P^{\otimes 2} \big( e^{\lambda_2 N_{0, m}} \big)\\ \label{eq:levrard2}
 &\to&
 e^{-\lambda_2}  \CstF
P^{\otimes 4} \left[
 e^{\Sigma_{0,m}}\right]   P^{\otimes 2} \big( e^{\lambda_2 N_{0, m}} \big) ,
\end{eqnarray}
by the local limit theorem (see e.g. \cite[Theorem 1.2.1, p.14]{Lawler}), with
\begin{equation}
\label{def:CstF}
\CstF = 2 \big(d/4\pi\big)^{d/2}.
\end{equation}
It follows that the sum in \eqref{eq:enfin} can be estimated, for $n \to \8$, by
\begin{eqnarray} \nn
n^{ {d-2}}  \sum_{k\geq n} {\rm of\ } \eqref{eq:enfin} &\sim&
 e^{-\lambda_2}  P^{\otimes 4} \left[
 e^{\Sigma_{0,m}}\right] \CstF  P^{\otimes 2} \big( e^{\lambda_2 N_{0, m}} \big)
n^{ {d-2}}  \sum_{k\geq n} \frac{1}{k^{d-1}} \\ \nn
&\to&
 e^{-\lambda_2}  \CstF P^{\otimes 4} \left[
 e^{\Sigma_{0,m}}\right]   P^{\otimes 2} \big( e^{\lambda_2 N_{0, m}} \big) \frac{1}{d-2} .
\end{eqnarray}
Collecting all this we conclude that the limit
$\lim_{n \to \8}T^<_2(m,n)$ exists, and further, that  the limits $\lim_{m \to \8} \lim_{n \to \8}T^<_2(m,n)$
and $\lim_{n \to \8}  T^<(n)$ exist and are equal.
Finally,
\begin{eqnarray} \nn
\lim_{n \to \8} \E s_n^4  &=& 2 \lim_{n \to \8} T^<(n) \\ \nn &\stackrel{\eqref{eq:enfin}}{=}&
2 \;\frac{\Cst \Zd \CstF(1\!-\!\pi_d)}{d-2} \; \Cstd e^{-\lambda_2} \times \E (W^4) \E (W^2) {\rm Var}(W)  \\ \label{eq:moment4} &=& \s_1^4 \E(W^4) .
\end{eqnarray}
This concludes the proof of \eqref{lem:instr4}.\medskip

\noindent

We then give the
{\bf proof of \eqref{lem:instr3}.}
We estimate the cross term
\begin{eqnarray}\nn
\E s_n^2 W_n^2
  & =& \Cstd   n^{\frac{d-2}{2}}  \sum_{ k\geq n}   P^{\otimes 4} \left[  e^{\Sigma_{0,n}+\lambda_2 N_{n,k}(S^{(1)},S^{(2)})} {\bf 1}_{S^{(1)}_{{k}} = S^{(2)}_{{k}}} \right]  \nonumber \\
&=&
 \Cstd   n^{\frac{d-2}{2}}  \sum_{ k\geq n}   P^{\otimes 4} \left[ e^{\Sigma_{0,n}} {\bf 1}_{S^{(1)}_{{k}} = S^{(2)}_{{k}}} P^{\otimes 2} \big[ e^{\lambda_2 N_{n,k}(S^{(1)},S^{(2)})} \big\vert \cS^{1,2}_{n,k}\big] \right]  , \nonumber
\end{eqnarray}
and we proceed as in \eqref{eq:levrard}, \eqref{eq:levrard2}.
On the event $\{S^{(1)}_{{k}} = S^{(2)}_{{k}}\} \bigcap \{S_{n-1}^{(i)}=s_{n-1}^{(i)}, i=1,2\}$, we have by time-reversal and Markov property,
\begin{eqnarray}
\nn
P^{\otimes 2} \big( e^{\lambda_2 N_{n, k}(S^{(1)}, S^{(2)})} \big \vert  \cS^{(1,2)}_{n,k}\big)
&=&
e^{-\lambda_2} P^{\otimes 2} \big( e^{\lambda_2 N_{0, k-n}} \big \vert  S^{(i)}_{k-n}=s^{(i)}_{n-1}, i=1,2\big)
\\
&\to& e^{-\lambda_2}  P^{\otimes 2} \big( e^{\lambda_2 N_{0, \8}} \big) , \nn
\end{eqnarray}
as $k \to \8$, and uniformly in $k \geq n$ as $n \to \8$,
$$
k^{\frac{d}{2}}  P^{\otimes 4} \left[  e^{\Sigma_{0,n}+\lambda_2 N_{n,k}(S^{(1)},S^{(2)})} {\bf 1}_{S^{(1)}_{{k}} = S^{(2)}_{{k}}} \right]
\to e^{-\lambda_2}  \E W^2 \CstF \E W^4 .
$$
Then,
\begin{eqnarray}\nn
\lim_{n \to \8} \E s_n^2 W_n^2
&=&
\Cstd e^{-\lambda_2}  \E W^2 \CstF \E W^4
 \lim_{n \to \8}
n^{\frac{d-2}{2}}  \sum_{ k\geq n} k^{-\frac{d}{2}}\\
\label{eq:fin3}
&=& \frac{2 \CstF}{d-2} \times \Cstd e^{-\lambda_2} \times \E W^4 \E W^2 ,\nn
\end{eqnarray}
which ends the proof of \eqref{lem:instr3}.     \qed

We now give the
\medskip

\noindent
{\bf Proof of Lemma \ref{lem:vargas}.}
For the first limit, we put $q=(d+\eps)/(d-1+\eps)$ and we estimate
\begin{eqnarray}\nn
 \E  T^=(n) &=&
 \Cstd^2 n^{ d-2}  \sum_{k \geq n}  P^{\otimes 4} \left[  e^{\Sigma_{0,k}}
 {\bf 1}_{S^{(1)}_{{k}} =  S^{(2)}_{{k}}}  {\bf 1}_{S^{(3)}_{{k}} =  S^{(4)}_{{k}}} \right]\\ \nn
 &\stackrel{\mbox{H\"older}}{\leq}&
\Cstd^2 n^{ d-2}  \sum_{k \geq n}  P^{\otimes 4} \left[  e^{(d+\eps)\Sigma_{0,k}} \right]^{1/(d+\eps)}
P^{\otimes 2} [S^{(1)}_{{k}} =  S^{(2)}_{{k}}]^{2/q}\\ \nn
&\leq& C   n^{ d-2} \sum_{k \geq n} k^{-d/q}
\end{eqnarray}
with $C \in (0,\infty) $ a constant, by the local limit theorem. Being of order $n^{-\eps/(d+\eps)}$, the  last term vanishes as $n\to \8$.
\medskip

We now prove the second limit, using arguments which are similar to the ones above.
Since the difference is non-negative, the norm $\E | T^<_1(n) - T^<_2(m,n) | $ is equal to
\begin{eqnarray}\nn
C \!\!&\!\!
 n^{ {d-2}}  \sum_{k\geq n} P^{\otimes 4} \left[
 e^{\Sigma_{0,m}+ \lambda_2 N_{k \! - \! m, k}(S^{(1)}, S^{(2)})}
 \left\{ e^{\Sigma_{m,k}-\lambda_2 N_{k \! - \! m,k}(S^{(1)}, S^{(2)})}-1 \right\}
  {\bf 1}_{S^{(1)}_{{k}} \!  = \!  S^{(2)}_{{k}}}
 G(S^{(3)}_k \! - \! S^{(4)}_k)\right]
 \\
 =&\!\!\!\!
 C
 n^{ {d-2}}  \sum_{k\geq n} P^{\otimes 4} \left[
 e^{\Sigma_{0,m}+ \lambda_2 N_{k \! - \! m, k}(S^{(1)}, S^{(2)})}
  {\bf 1}_{S^{(1)}_{{k}} \!  = \!  S^{(2)}_{{k}}}
 G(S^{(3)}_k \! - \! S^{(4)}_k)\right.\qquad \qquad \qquad \qquad \nn\\ \label{eq:norm}
 & \qquad \qquad\qquad
 \left.
 P^{\otimes 4}
 \left\{ e^{\Sigma_{m,k}-\lambda_2 N_{k \! - \! m,k}(S^{(1)}, S^{(2)})}-1
 \big \vert
 S^{(i)}_t, 1\leq i\leq 4, t=m, k
 \right\}
 \right]  
\end{eqnarray}
with $C=  \frac{{\rm Var}(W)}{G(0)}  \Cstd$, by conditioning on the paths at times $t=1,\ldots m $ and $t=k$. The event
$$B_{m,k}= \big\{|S^{(i)}_t-S^{(j)}_t| \geq m^{1/4}; 1\leq i \neq j\leq 4, t=m, {\rm\ and\ } t=k \ {\rm except \ for \ } (i,j)=(1,2) \big\}$$ has probability larger than $1-c m^{-d/4}$ for some constant $c\in(0,\infty)$, by the local limit theorem.
On this event, by transience, the intersections between the paths between times $m$ and $k-m$ essentially come, when $m$ and $k$ are large,
from those of $S^{(1)}$ and $S^{(2)}$ between times $k-m$ and $m$. Precisely,
\begin{equation}
\lim_{m \to \8} \limsup_{k \to \8}
\sup_{B_{m,k}}  P^{\otimes 4}  \left\{ \Sigma_{m,k} \neq \lambda_2 N_{m \! - \! k,k}(S^{(1)}, S^{(2)})
 \big \vert  S^{(i)}_t, 1\! \leq\! i \! \leq \!4, t\!=\!m, k \right\} =0 \nn
\end{equation}
(note that the event in the previous line means that there are no intersections of $S^{(i)}$ and $S^{(j)}$ between times $m$ and $k$ for $i\geq 3$ or $j \geq 3$,
and  no intersections of $S^{(1)}$ and $S^{(2)}$ between times $m$ and $k-m$). Under our integrability condition \eqref{eq:int1}, this is enough to imply that
\begin{equation}\nn
\lim_{m \to \8} \limsup_{k \to \8}
\sup_{B_{m,k}}  P^{\otimes 4}  \left\{ e^{\Sigma_{m,k}-\lambda_2 N_{m \! - \! k,k}(S^{(1)}, S^{(2)})}-1
 \big \vert  S^{(i)}_t, 1\! \leq\! i \! \leq \!4, t\!=\!m, k \right\} =0.
\end{equation}
Plugging this in \eqref{eq:norm} and using the fact that $ n^{ {d-2}}  \sum_{k\geq n} P^{\otimes 4} [
 e^{\Sigma_{0,m}+ \lambda_2 N_{k \! - \! m, k}(S^{(1)}, S^{(2)})}
  {\bf 1}_{S^{(1)}_{{k}} \!  = \!  S^{(2)}_{{k}}}]$ has a finite limit,
we conclude that the right-hand side of \eqref{eq:norm} vanishes as $n\to \8$. \qed

\bigskip
\noindent
{\bf Proof of Lemma \ref{lem:rate-conv-Dk}}.
With the notation introduced in the beginning of this section, we have
$D_{k+1} = W_{k+1} - W_k = P\big[ e^{h_k(S)} (e^{ \beta \eta (k,S_k) - \lambda (\beta)} -1)\big]$,
$$ D^4_{k+1}   = P^{\otimes 4} \left[ e^{\sum_{i=1}^4 h_k(S^{(i)} )}  \prod_{i=1}^4 (e^{ \beta \eta (k,S_k^{(i)}) - \lambda (\beta)} -1) \right].
                                      $$
Using Fubini's theorem and independence,  we obtain
\begin{eqnarray}
\E D^4_{k+1}
   &=&   P^{\otimes 4}  \left[ \E e^{\sum_{i=1}^4 h_k(S^{(i)} )}  \E \prod_{i=1}^4 (e^{ \beta \eta (k,S_k^{(i)}) - \lambda (\beta)} -1) \right] \nonumber  \\
    &=& P^{\otimes 4}  e^{\Sigma_{0,k}}  (\gamma_4 {\bf 1}_{ (S_k^{(1)},  S_k^{(2)}, S_k^{(3)}, S_k^{(4)}) \in E_4 }
      + \gamma_2   {\bf 1}_{ (S_k^{(1)},  S_k^{(2)}, S_k^{(3)}, S_k^{(4)}) \in E_{2,2} } ),
\end{eqnarray}
where $ \gamma_4 = \E (e^{ \beta \eta (0,0)  - \lambda (\beta)} -1)^4$ and  $ \gamma_2 = (\E (e^{ \beta \eta (0,0) - \lambda (\beta)} -1)^2)^2$. Notice that   $(S_k^{(1)},  S_k^{(2)}, S_k^{(3)}, S_k^{(4)}) \in E_{2,2}$ if and only if one of the following cases occurs: (a) $ S_k^{(1)} = S_k^{(2)} \neq S_k^{(3)}= S_k^{(4)}$,
 (b) $ S_k^{(1)} = S_k^{(3)} \neq S_k^{(2)}= S_k^{(4)}$, (c) $ S_k^{(1)} = S_k^{(4)} \neq S_k^{(2)}= S_k^{(3)}$. Therefore by symmetry, we obtain
$$ \E D^4_{k+1}   \leq  (\gamma_4 + 3\gamma_2) P^{\otimes 4}e^{\Sigma_{0,k}}  {\bf 1}_{S_k^{(1)} = S_k^{(2)}}
   {\bf 1}_{S_k^{(3)} = S_k^{(4)}}.   $$
Hence, using H\"older's inequality, for $p, q>1$ with $1/p + 1/q =1$, we have
$$ \E D^4_{k+1}   \leq  (\gamma_4 + 3\gamma_2)
 (P^{\otimes 4}e^{ p \Sigma_{0,k}})^{1/p}  (P^{\otimes 2} (S_k^{(1)} = S_k^{(2)} ))^{2/q}.   $$
By the local limit theorem, $P^{\otimes 2} (S_k^{(1)} = S_k^{(2)} )= O (k^{-d/2})$. Thus taking
$ |\beta| >0$ small enough such that $P^{\otimes 4}e^{ p \Sigma_{0,k}} < \infty$, we see that
Eq. \eqref{eq-rate-conv-Dk} holds. This ends the proof of Lemma \ref{lem:rate-conv-Dk}. \qed

\bigskip

\noindent
{\bf Acknowledgements.} The work  is supported in parts by  CNRS (UMR 7599
Probabilit{\'e}s et Mod{\`e}les
Al{\'e}atoires et UMR 6205 Laboratoire de Math\'ematique de Bretagne Atlantique), the National Natural Science Foundation of China (Grants no. 11571052 and no. 11401590) and
the U.S. National Science Foundation (Grant no. NSF PHY11-25915).

{\small


}


\end{document}